\documentclass[runningheads]{llncs}
\usepackage[T1]{fontenc}
\usepackage{graphicx}
\usepackage{amsmath,amscd,amssymb}
\usepackage{listings}

\DeclareMathOperator{\proj}{proj}
\DeclareMathOperator{\fr}{fr}
\DeclareMathOperator{\Frontier}{Frontier}
\DeclareMathOperator{\cl}{cl}

\newcommand{\lex}{<_{\rm{lex}}}
\newcommand{\lexeq}{\le_{\rm{lex}}}
\newcommand{\R}{{\mathbb R}}

\newtheorem{algorithm}[theorem]{Algorithm}

\begin{document}
\title{Cylindrical Algebraic Decompositions with Frontier Condition\thanks{The author was supported by an EPSRC Doctoral Training grant}}
%
%
\author{Hollie Baker} 
%
\authorrunning{H. Baker}
%
\institute{University of Bath, Bath, BA2 7JY, UK.
\email{hb581@bath.ac.uk}}
%
\maketitle              
\begin{abstract}
A Cylindrical Algebraic Decomposition (CAD) is a decomposition of $\R^n$ into a finite collection of semialgebraic cells. A CAD satisfies the ``frontier condition'' if, for every cell $C$, there is a collection of cells of the decomposition whose union is the closure of $C$. This property is referred to in other literature as ``closure finiteness'' or ``boundary coherence''.
This paper proves the {\em existence} of, and presents an algorithm to construct, a CAD satisfying the frontier condition without a preliminary change of coordinates, e.g., in the potential presence of blow-ups.
The algorithm has {\em elementary} (in the sense
of L. Kalmar) complexity. This also provides an upper bound on the number of cells in a CAD with this property. The frontier condition can be useful in computing topological properties of semialgebraic sets defined by first-order formulas, in solving motion planning problems and in triangulations
of definable monotone families.
The algorithm presented takes a novel approach in that it uses a recursion on the lexicographical index of cells in the initial decomposition.

\keywords{cylindrical algebraic decomposition \and complexity \and frontier condition \and closure finite \and boundary coherent}
\end{abstract}
%
%

\section{Introduction}

We prove that a Cylindrical Algebraic Decomposition (CAD) of $\R^n$, compatible with a semialgebraic set $S \subset \R^n$ and satisfying the frontier condition, can be constructed in any dimension and without a preliminary change of coordinates even in the presence of blow-ups. A CAD is compatible with a set $S$ if $S$ is the union of some cells of the CAD. Frontier condition means that the closure of every cell is a union of some cells in the decomposition.
To our knowledge, this is the first proof that such a CAD exists without the change of coordinates.
We present the proof in the form of an algorithm which constructs a CAD compatible with a semialgebraic set $S \subset \R^n$ and satisfying the frontier condition.
An upper bound on complexity is obtained. This is also an upper bound on the number of cells, number of polynomials and degree of polynomials in the CAD. The algorithm has {\em elementary} complexity (in the sense of L. Kalmar, see e.g.,
\cite[$\S 57$]{kleene1952}).

The frontier condition is useful in computing topological properties of semialgebraic sets defined by first-order Boolean formulas. For example, these decompositions can be viewed as closure-finite weak cell complexes (CW-complexes) and their homologies can be computed. In addition, cell adjacencies can easily be computed and motion-planning problems, for instance the {\em piano-mover's problem} presented by Schwartz and Sharir in \cite{pianomovers1983}, can be solved, leading to applications in robotics.

It was shown by Schwartz and Sharir in \cite{pianomovers1983} that it is relatively easy to satisfy the frontier condition after a generic linear rotation of coordinates in $\R^n$. However, such a rotation of coordinates may be undesirable, or even impossible, in some applications.

An important application relates to the work of Basu, Gabrielov and Vorobjov in \cite{bgv15}, who present an algorithm for computing the triangulation of a definable monotone family. This algorithm requires a CAD of the family such that every cell is topologically regular and satisfies the frontier condition. Since the definable monotone family depends on the initial coordinate ordering, a rotation of coordinates is not allowed. Davenport, Locatelli and Sankaran \cite{jhd20} prove that if a CAD is {\em well-based} then it satisfies these properties.
A CAD is well-based if it contains no {\em bad cells}: cells above which one of the polynomials defining the CAD vanishes identically. These occur frequently in definable monotone families. Therefore, an algorithm for computing CAD with frontier condition, compatible with sets containing blow-ups and without a change of coordinates, is needed. Such a CAD algorithm is presented by the authors \cite{bgv15} for sets of dimension not greater than two. The algorithm presented in this paper provides a step towards extending their result to decompositions compatible with sets (and potentially triangulations of definable monotone families) of arbitrary dimension.

In addition, the result can be applied to other categories, e.g., semialgebraic sets defined by {\em fewnomials} (see Section~\ref{sec:pfaffian}), whose structure is destroyed by a change of coordinates.

The result can also be easily extended to semialgebraic sets defined by first-order Boolean formulas with quantifiers.
Indeed, by Theorem~\ref{th:proj}, the projection of a CAD with frontier condition to any dimension also satisfies the frontier condition, and quantifier elimination requires considering projections of the CAD.

It is important to note that D. Lazard \cite{lazard10} presented an algorithm for obtaining a CAD in $\R^n$ for $n \le 3$ satisfying the frontier condition without a change of coordinates. Lazard takes an algebraic approach and the complexity of the algorithm presented is the same as that of classical CAD. Our result takes a different approach, relying on a recursion on the lexicographical order of cell indices, and works for CAD of arbitrary dimension.

The complexity results obtained in this paper are of complexity-theoretical nature, and are not claimed to be
useful in practical computations.
However, we emphasise that this algorithm has {\em elementary} complexity, being triply exponential in the number of variables.
\bigskip

We now introduce, closely following \cite{bgv15}, some notations and definitions needed to formulate the main result of the paper.

For a semialgebraic set $X \subset \R^n$, let $\cl (X)$ denote the closure of $X$ in Euclidean topology and
$\fr (X):= \cl (X) \setminus X$ be the {\em frontier} of $X$.
Let $\proj_{\R^k}:\> \R^n \to \R^k$, for some $1 \le k \le n$, be the projection map from $\R^n$ to $\R^k$.

\begin{definition}[cylindrical cells]
\label{def:cells}
Let $n \ge 1$ and $(i_1, \ldots, i_n ) \in \{0, 1\}^n$.
A cylindrical $(i_1 ,\ldots, i_n )$-cell is a semialgebraic set in $\R^n$ obtained by induction
on $n$ as follows.
\begin{itemize}
\item $n=1$: A $(0)$-cell is a single point $x \in \R$, a $(1)$-cell is an open interval in $\R$: either $(x, y)$, $(-\infty, y)$, $(x, \infty)$,
or $(-\infty, \infty)$ ($x,y \in \R$).
\item $n>1$: Suppose that $(i_1, \ldots, i_{n-1})$-cells, are deﬁned.
An $(i_1, \ldots, i_{n-1}, 0)$-cell (or a section cell) is the graph in $\R^n$ of a polynomial $f : C' \to \R$,
where $C'$ is an $(i_1, \ldots, i_{n-1})$-cell.
An $(i_1 , \ldots, i_{n-1}, 1)$-cell (or a sector cell) is a subset of $C' \times \R$: either $\{(x, t) \in C' \times \R \mid f (x) < t < g(x)\}$, $\{(x, t) \in C' \times \R \mid f (x) < t\}$, $\{(x, t) \in C' \times \R \mid t < g(x)\}$, or $C' \times \R$, where $C'$ is an $(i_1, \ldots, i_{n-1})$-cell and $f, g:\> C' \to \R$ are polynomials such that $f (x) < g(x)$ for all $x \in C'$.
\end{itemize}
\end{definition}

\begin{definition}
\label{def:sector-top-bottom}
Let $C$ be a sector cell with polynomials $f,g$ as in Definition~\ref{def:cells}. The graph of $f$ is called the {\em bottom} of $C$, and the graph of $g$ is called the {\em top} of $C$.
\end{definition}

\begin{definition}
\label{def:section-top-bottom}
Let $C$ be a $(i_1, \ldots, i_{n-1}, 0)$- (section-)cell, and $k$ be the largest among $\{1, \ldots, n - 1\}$ such that $i_k = 1$.
Then $C$ is the graph of a map $\mathbf{f} :\>C' \to \R^{n-k}$, where $C'$ is a $(i_1, \ldots, i_k)$- (sector-)cell.
The pre-image of the bottom of $C'$ by $\proj_{\R^k}\vert_{\cl(C)}$ is called the bottom of $C$,
and the pre-image of the top of $C'$ by $\proj_{\R^k}\vert_{\cl(C)}$ is called the top of $C$.
\end{definition}

Informally, the bottom (resp. top) of a section cell consists of the points of the closure of $C$ ``above'' the bottom (resp. top) of the sector cell $C'$.
In some literature, section cells are called {\em graphs} and sector cells are called {\em bands}.

It is clear from the definition that the top and bottom of a sector cell, if they exist, are always cylindrical section cells. However, this is not
always the case for a section cell. We now present some examples of $2$-dimensional section cells in $\R^3$ to illustrate. In Example~\ref{example:nice-cell}, a cylindrical cell whose top and bottom are cylindrical cells is presented. Example~\ref{example:bottom-vanish} shows that the top or bottom of a cell may fail to exist and Example~\ref{example:top-bottom-not-cylindrical} shows that, even if it exists, the top or bottom of a cell may not be a cylindrical cell.

\begin{example}\cite[Example 4.3]{bgv13sets} \label{example:nice-cell}
Let $n=3$, $C' := \{ (x,y) \mid 0 < x < 1, y > 0, x + y < 1 \}$, $f(x,y) = x^2 + y^2$, and $C$ be the graph of $f\vert_{C'}$.
Observe that the bottom ${C'}_B = \{ (x,y) \mid 0 < x < 1, y=0 \}$ and top ${C'}_T = \{ (x,y) \mid 0 < x < 1, x+y=1 \}$ of $C'$ are 1-dimensional section cells. Consider
$$C_B = \left\{ (x,y,z) \mid (x,y)\in {C'_B}, z = x^2 + 0 \right\}$$ and $$C_T = \left\{ (x,y,z) \mid (x,y)\in {C'_T}, z = x^2 + y^2\right\}.$$ Both $C_B$ and $C_T$ are graphs of polynomials, hence they are 1-dimensional section cells ($(1,0,0)$-cells).
\end{example}

Sometimes the top and bottom may not exist. This could happen if ${C'}_T$ or ${C'}_B$ fail to exist, i.e., $C'$ is not bounded from below, or from above. The top (resp. bottom) of a section cell may also be empty if the graph of $f : C' \to \R$ tends to infinity as it approaches every point in ${C'}_T$ (resp. ${C'}_B$), as illustrated in the following example.

\begin{example} \label{example:bottom-vanish}
Let $n=3$, $C' := \{ (x,y) \mid x > 0, y > 0 \}$, $g(x,y) = 1/x$ and $C$ be the graph of $g\vert_{C'}$.
Note that $C = \{ (x,y,z) \mid x > 0, y > 0, x z - 1 = 0 \}$, so $C$ a semialgebraic set.
On the other hand, $xz - 1 = 0$ is undefined at every point in ${C'}_B = \{ (x,0) \mid x>0 \}$ and $\{ xz - 1 = 0 \}$ tends to infinity as $y\to 0$ in $C$.
Hence, $C_B$ does not exist. Neither does $C_T$, since ${C'}_T$ does not exist as $C'$ is not bounded from above.
\end{example}

Finally, let us consider an example where $C_B$ exists but is not a cylindrical cell.
This can happen when ${C'}_B$ contains a blow-up point -- a point at which some polynomial vanishes identically.

\begin{example} \label{example:top-bottom-not-cylindrical}
\cite[Example 3.2]{bgv15}
Let $n=3$, $C' = \{ (x,y) \mid -1 < x < 1, \vert x\vert < y < 1 \}$.
Let $\varphi(x,y) = \vert x/y\vert$, and let $C$ be the graph of $\varphi\vert_{C'}$.
The bottom of $C$ is not the graph of a continuous function. In particular, due to the blow-up point of $\varphi$ at the origin,
$$C_B = \{ (x,y,z) \mid -1 < x < 1, x \ne 0, y = \vert x\vert, z = 1 \} \cup \{(0,0,z) \mid 0\le z \le 1\}.$$
Note that $C_B$ is not the image of $\varphi$ because $\varphi$ vanishes identically at the origin, but $\fr(C) \cap \{ (0,0,z) \mid z \in \R \} = \{ (0,0,z) \mid 0 \le z \le 1 \}$.
\end{example}

\begin{figure}
\begin{center}
\includegraphics[width=0.5\textwidth]{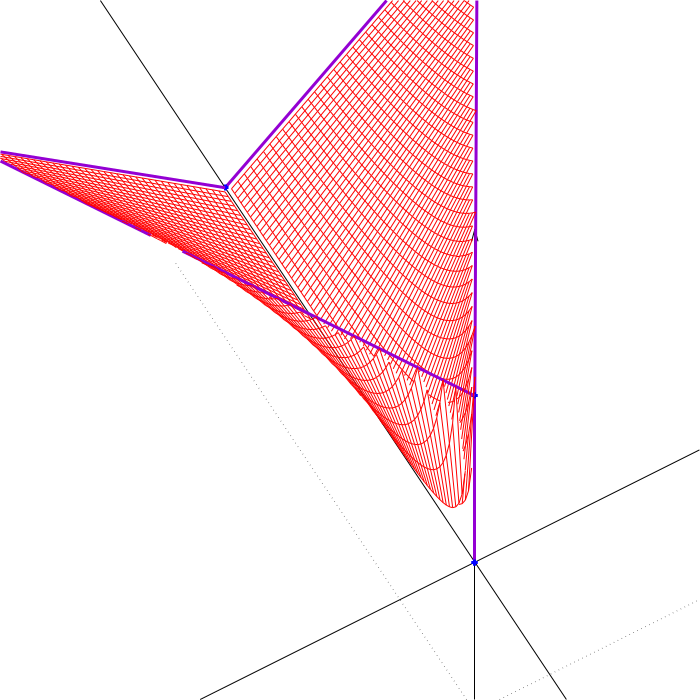}
\end{center}
\caption{shows $C$ and $\fr(C)$ from Example~\ref{example:top-bottom-not-cylindrical}, demonstrating that $C_B$ is not the graph of a continuous function.}
\end{figure}

\begin{definition}[CAD]
A cylindrical cell decomposition of $\R^n$ is a finite partition of $\R^n$ into
cylindrical cells defined by induction on $n$ as follows.
\begin{itemize}
\item When $n=0$, the cylindrical cell decomposition of $\R^n$ consists of a single point.
\item Let $n > 0$.
For a partition $\mathcal D$ of $\R^n$ into cylindrical cells, let ${\mathcal D}'$
be the set of all cells $C' \subset \R^{n-1}$ such that $C' = \proj_{\R^{n-1}} (C)$ for some cell $C$ of $\mathcal D$.
Then $\mathcal D$ is a cylindrical cell decomposition of
$\R^n$ if ${\mathcal D}'$ is a cylindrical cell decomposition of $\R^{n-1}$.
In this case, we call ${\mathcal D}'$ the cylindrical cell decomposition of $\R^{n-1}$ induced by $\mathcal D$.
\end{itemize}
\end{definition}

\begin{definition}
\begin{enumerate}
\item
A cylindrical cell decomposition $\mathcal D$ of $\R^n$ is compatible with a semialgebraic
set $X \subset \R^n$ if, for every cell $C$ of $\mathcal D$, either $C \subset X$ or $C \cap X = \emptyset$.
\item
A cylindrical cell decomposition ${\mathcal D}'$ of $\R^n$ is a {\em refinement} of a decomposition $\mathcal D$
of $\R^n$ if ${\mathcal D}'$ is compatible with every cell of $\mathcal D$.
\end{enumerate}
\end{definition}

\begin{definition} \cite[Definition 3.12]{bgv15}
A decomposition $\mathcal{D}$ satisfies the frontier condition if for each cylindrical cell $C$, its frontier $\cl(C) \setminus C$ is a union of cells of $\mathcal{D}$.
\end{definition}

As discussed above, in relation to Example~\ref{example:top-bottom-not-cylindrical}, the frontier of a cylindrical cell
can be geometrically non-trivial, e.g., it may contain a blow-up.

We now make precise how a semialgebraic set $S\subset \R^n$ will be presented so that we can state our main theorem.

\begin{definition} \label{def:qff}
A quantifier-free Boolean formula is obtained by the following rules:
\begin{enumerate}
\item \label{item:qff-1} If $f \in \R[x_1\ldots,x_n]$ then $f=0$ and $f>0$ are quantifier-free Boolean formulas,
\item \label{item:qff-2} if $F$ and $G$ are quantifier-free Boolean formulas then $F \land G$, $F \lor G$ and $\neg F$ are quantifier-free Boolean formulas.
\end{enumerate}
\end{definition}

\begin{definition}
A first-order Boolean formula is obtained by rules \ref{item:qff-1} and \ref{item:qff-2} from Definition~\ref{def:qff} and
\begin{enumerate}
  \setcounter{enumi}{2}
  \item If $F$ is a first-order Boolean formula and $y$ is a variable ranging over $\R$, then $\exists y\ F$ and $\forall y\ F$ are first-order Boolean formulas.
\end{enumerate}
\end{definition}

By the Tarski-Seidenberg theorem, every quantifier-free Boolean formula represents a semialgebraic set, and every semialgebraic set can be represented by a quantifier-free Boolean formula. We will use quantifier-free Boolean formulas as the means of representing semialgebraic sets in the following satements.

\begin{theorem}\label{th:main}
Let $S \subset \R^n$ be a semialgebraic set defined by a quantifier-free Boolean formula $F$ with $s$ different polynomials of maximum degree $d$ in $\R[x_1,\ldots,x_n]$.
There is an algorithm, taking $F$ as input, which outputs a cylindrical decomposition $\mathcal D$ of $\R^n$ compatible with $S$ and
satisfying the frontier condition.
The complexity of this algorithm is $(sd)^{O(1)^{n2^n}}$.
This is also an upper bound on the number of cells in $\mathcal D$, number of polynomials defining cells and their degrees.
\end{theorem}

\begin{remark}
Algorithm in Theorem~\ref{th:main} is understood as a Blum-Shub-Smale (BSS) real numbers machine \cite{blum1998}.
A similar statement is also true for the Turing machine model, in which case the complexity bound
depends, in addition, polynomially on maximal bit-size of coefficients (cf. \cite{collins1975}).
\end{remark}

We will also prove the following.

\begin{theorem}
\label{th:proj}
Let $\mathcal{D}$ be a CAD of $\R^n$ compatible with a semialgebraic set $S\subset \R^n$ and satisfying the frontier condition.
Then the decomposition induced by $\mathcal{D}$ on $\R^{k}$ for all $1\le k < n$ is a cylindrical decomposition of $\R^k$ compatible with $\proj_{\R^k}(S)$ and satisfying the frontier condition.
\end{theorem}

This result allows us to apply Theorem~\ref{th:main} to semialgebraic sets represented by first-order Boolean formulas (with quantifiers) and obtain a decomposition compatible with $S$ and satisfying the frontier condition.

\section{Subroutines}\label{sec:frontier}

\subsection{Classical CAD}
The following statement was proved in \cite{collins1975} (and almost simultaneously and independently in \cite{wuthrich2005}). Here we formulate a slightly improved version.

\begin{proposition}\label{prop:collins}\cite[Algorithm~11.2]{bpr2006}
Let $f_1,\ldots,f_s$ be polynomials in $\R [x_1, \ldots ,x_n]$ having maximum degree $d$.
There is an algorithm, taking $\{f_1,\ldots,f_s\}$ as input, which produces a cylindrical decomposition $\mathcal E$ of $\R^n$ such that every cell in $\mathcal{E}$ has constant sign on every polynomial $f_1,\ldots,f_s$.
The complexity of the algorithm is $(sd)^{O(1)^{n}}$.
This is also an upper bound on the number of cells in $\mathcal E$, number of polynomials defining cells and their degrees.
\end{proposition}

We now describe two variants of the CAD algorithm.

\begin{proposition}\label{prop:collins-sets}
Let $S_1,\ldots,S_k$ be a finite collection of semialgebraic subsets of $\R^n$ defined by quantifier-free Boolean formulas $F_1,\ldots,F_k$ respectively. Together, these formulas include $s$ different polynomials in $\R[x_1,\ldots,x_n]$, having maximum degree $d$.
There is an algorithm, taking $\{F_1,\ldots,F_k\}$ as input, which produces a cylindrical decomposition $\mathcal E$ of $\R^n$ compatible with each $S_i, 1 \le i \le k$.
Complexity, number of cells, number of polynomials and degrees are the same as in Proposition~\ref{prop:collins}.
\end{proposition}

\begin{proof}
Let $f_1,\ldots,f_s \in \R[x_1,\ldots,x_n]$ be the polynomials from formulas $F_1,\ldots,F_k$ and let $\mathcal{E}'$ be the CAD of $\R^n$ produced by Proposition~\ref{prop:collins} with $f_1,\ldots,f_s$ as input. Since $F_i$ for each $1\le i\le k$ provides a set of sign conditions on polynomials $f_1,\ldots,f_s$, $\mathcal{E}'$ is a refinement of the CAD produced by this algorithm.

Truth values of $F_1,\ldots,F_k$ on ecah cell can easily be computed from the CAD produced by Proposition~\ref{prop:collins} by considering the signs of $f_1,\ldots,f_s$ on ecah cell, which we computed in the construction of $\mathcal{E}'$.
Observe that this does not change the asymptotic complexity or bounds in Propoosition~\ref{prop:collins}.
\end{proof}

The construction described in the above proof is very naive. A description of this algorithm is presented in, e.g., \cite{collins1991}, and an alternative approach is presented in \cite{bradford2014}.

\begin{proposition}\label{prop:collins-set}
Let $S\subset \R^n$ be a semialgebraic set defined by a quantifier-free Boolean formula $F$ containing $s$ different polynomials in $\R[x_1,\ldots,x_n]$, having maximum degree $d$.
There is an algorithm, taking $F$ as input, which produces a cylindrical decomposition $\mathcal E$ of $\R^n$ compatible with $S$.
Complexity, number of cells, number of polynomials and degrees are the same as in Proposition~\ref{prop:collins}.
\end{proposition}

\begin{proof}
Observe that this is the same as the algorithm as in Proposition~\ref{prop:collins-sets} taking $\{F\}$ as input.
\end{proof}

\subsection{Computing the Frontier}

Now we describe an algorithm, having singly-exponential complexity, for finding the frontier, $\fr(S)=\cl (S) \setminus S$,
of a semialgebraic set $S \subset \R^n$.

\begin{lemma}\label{le:frontier}
Let $S \subset \R^n$ be a semialgebraic set defined by a quantifier-free Boolean formula $F$
containing $s$ different polynomials in $\R [x_1, \ldots ,x_n]$ having maximum degree $d$.
There is an algorithm, taking $F$ as input, which represents the semialgebraic set $\fr (S)$ by a quantifier-free Boolean formula $F'$ with complexity $(sd)^{O(n^2)}$.
This is also an upper bound on the number of polynomials in $F'$ and their degrees.
\end{lemma}

\begin{proof}
Observe that $\fr(S)$ can be represented by a first-order Boolean formula
$$
\fr(S) = \left\{ \mathbf{x} \in (\R^n \setminus S) \mid \forall \varepsilon >0\> \exists \mathbf{y} \in S\> (\Vert \mathbf{x} - \mathbf{y} \Vert < \varepsilon ) \right\}.
$$
Using singly-exponential quantifier elimination algorithm \cite[Algorithm~14.21]{bpr2006}, we represent $\fr (S)$
as a quantifier-free Boolean formula $F'$ with the bounds required in the lemma.
\end{proof}

\section{Obtaining a CAD with Frontier Condition}

Recall that, according to Definition~\ref{def:cells}, to each cylindrical cell in a CAD of $\R^n$ a (multi-)index
$(i_1, \ldots ,i_n) \in \{ 0,1 \}^n$ is assigned.
Introduce a lexicographic order $\lex$ on the set of all cell indices as follows.
For any two indices, $M:=(i_1,\ldots,i_n)$ and $N:=(j_1,\ldots,j_n)$, we set $M \lex N$
iff for the maximal $k,\ 0 \le k <n$, such that $i_1=j_1, \ldots, i_k=j_k$ we have $i_{k+1} < j_{k+1}$.
If $M \lex N$ or $M = N$, we write $M \lexeq N$.

\begin{lemma} \label{lem:fr-lex-less}
Let $\mathcal{D}$ be a CAD of $\R^n$ and let $C$ be a cell of $\mathcal{D}$ with index $M$.
Then $\fr(C)$ is contained in a union of cells of $\mathcal{D}$ with indices lexicographically less than $M$.
\end{lemma}

\begin{proof}
Proceed by induction on $n$.

When $n=1$, the cell $C$ either has index $M=(0)$ or $M=(1)$.
If $M=(0)$ then $C$ is a single point and its frontier is empty.
If $M=(1)$, $C$ is an open interval. By the definition of CAD, endpoints of $C$ are $(0)$-cells.

Now let $n>1$ and let $C$ have index $M=(i_1,\ldots,i_n)$.
The projection $C'=\proj_{\R^{n-1}}(C)$ is a cell in the induced decomposition $\mathcal{D'}$ of $\R^{n-1}$
with index $(i_1,\ldots,i_{n-1})$.
By the induction hypothesis, $\fr(C')$ is contained in a union of cells of $\mathcal{D'}$ with indices
$(j_1,\ldots,j_{n-1}) \lex (i_1,\ldots,i_{n-1})$.

If $i_n=0$, then $C$ is a section cell.
Its frontier $\fr (C)$ is contained in $\fr (C') \times \R$ and, therefore, in a union of
cells of $\mathcal D$ with indices $(j_1,\ldots,j_{n-1}, j_n)$ for some $j_n \in \{ 0,1 \}$. We have shown that $(j_1,\ldots,j_{n-1}, j_n) \lex M$.

If $i_n=1$, then $C$ is a sector cell.
In this case, $\fr (C)$ is contained in $(\fr (C') \times \R) \cup C_T \cup C_B$, hence in the
union of cells of $\mathcal D$ with indices $(j_1,\ldots,j_{n-1}, j_n)$ for some $j_n$
and two cells, $C_T$ and $C_B$, having the same index $(i_1,\ldots,i_{n-1},0)$.
All of these indices are lexicographically less than $M$.
\end{proof}

In further proofs we will need the following technical statement.

\begin{lemma}\label{lem:section-split}
Let $f: X \to \R$ be a continuous polynomial function on an open set
$X \subset \R^n$ having the graph $F \subset \R^{n+1}$.
Let $G$ be a semialgebraic subset of $F$, which is closed in $F$, and having dimension less than $\dim (F)=n$.
Then $\fr_{F}(F \setminus G) = G$, where $\fr_F$ denotes the frontier in $F$.
\end{lemma}

\begin{proof}
Let $Y = \proj_{\R^{n-1}}(G)$ and observe that $X = \proj_{\R^{n-1}}(F)$.
Let us first prove that $\fr_X(X \setminus Y) = Y$.
By the definition of frontier, we have
$$\fr_X(X \setminus Y)=\cl_X(X \setminus Y) \setminus (X \setminus Y),$$
where $\cl_X$ denotes the closure in $X$.
$\cl_X(X \setminus Y)=X$. Indeed, $X\setminus Y$ is open and $n$-dimensional since $Y$ is closed in $X$ and $\dim(Y)<n$.
Assume that $\cl_X(X \setminus Y) \ne X$. Then there is a nonempty subset $Z\subset X$ such that $Z \not \subset \cl_X(X \setminus Y)$. Since $X \setminus Y \subset \cl_X(X \setminus Y)$, $Z$ is necessarily a subset of $Y$ and $Z$ has dimension $<n$. This implies that $\cl_X(X\setminus Y)$ is not, which is a contradiction.
It follows that $\fr_X(X \setminus Y)=X \setminus (X \setminus Y)= X \cap Y = Y$.

Finally, since $f$ is a continuous function, and therefore $Y$ and $X$ are homeomorphic to $G$ and $F$
respectively, we obtain the required formula $\fr_{F}(F \setminus G) = G$.
\end{proof}

We now prove that the decomposition induced on $\R^k, 1 \le k \le n-1$ by any CAD satisfying the frontier condition also satisfies the frontier condition.

\begin{proof}[of Theorem~\ref{th:proj}.]
Let $\mathcal{D}$ be a CAD compatible with $S$ and satisfying the frontier condition.
By the definition of CAD, $\proj_{\R^k0}(S)$ is a union of cells in the decomposition induced by $\mathcal{D}$ on $\R^k$ for all $1\le k \le n-1$.

We now prove that if $\mathcal{D}$ satisfies the frontier condition then all decompositions induced by $\mathcal{D}$ on $\R^k$ for all $1 \le k \le n-1$ satisfy the frontier condition.
Consider $\Gamma$, the union of cells of $\mathcal{D}$ contained in the $k$-dimensional cylinder $\proj_{\R^{n-1}}(C) \times \R$ where $C$ is a cell of $\mathcal{D}$.
Observe that $\bigcup_{C \in \Gamma} C = \proj_{\R^{n-1}}(C) \times \R$.
Consider $$
\fr(\Gamma) := \bigcup_{C \in \Gamma} \fr(C) \setminus \bigcup_{C \in \Gamma} C . 
$$
Observe that $\fr(\Gamma)$ is a union of cells of $\mathcal{D}$ since every cell $C$ in $\Gamma$ satisfies the frontier condition. 
We claim that $\fr(\Gamma)$ is the union of some similarly defined cylinders of dimension less than $k$.
Proceed by induction on $k$.

The base case, $k=1$, is simple since $\fr(\Gamma) = \emptyset$.

Now suppose that $k>1$. Each sector cell in $\Gamma$ has dimension $k$ and each section cell in $\Gamma$ has dimension $k-1$, so the dimension of the union of frontiers of all cells in $\Gamma$ is at most $k-1$.
Hence $\dim(\fr(\Gamma))$ is at most $k-1$.
By the cylindrical property of CAD, the projection of any two cells contained in $\fr(\Gamma)$ is either disjoint or coincides. Therefore, $\fr(\Gamma)$ is a union of cylinders $\Gamma' := \proj_{\R^{n-1}}(B) \times \R$ where $B$ is a cell of $\mathcal{D}$ contained in $\fr(\Gamma)$. Observe tthat each cylinder $\Gamma'$ has dimension less than $k$.
By the induction hypothesis, $\fr(\Gamma')$ is a union of cylinders with dimension less than $\dim(\Gamma')$.

By the definition of CAD, each of these cylinders $\Gamma'$ projects onto a single cell in the decomposition $\mathcal{D'}$ induced by $\mathcal{D}$ on $\R^{n-1}$. Hence $\fr(C')$, where $C' := \proj_{\R^{n-1}}(C)$ for all $C \in \Gamma$ is a union of cells $B' := \proj_{\R^{n-1}}(B)$ of $\mathcal{D'}$ for all $B \in \Gamma'$.
It follows that $\mathcal{D'}$ satisfies the frontier condition.

The proof that every decomposition induced by $\mathcal{D}$ on $\R^{k}, 1 \le k \le n-1$ satisfies the frontier condition is completed by iterating the projection operation.
\end{proof}

We are now ready to prove the main result. The worked example, presented in Section~\ref{sec:worked-example} is referenced to help illustrate the procedure.

\begin{proof}[of Theorem~\ref{th:main}]
We begin with a mathematical description of the algorithm.
\medskip

\noindent{\bf Algorithm.}

{\bf Input:} $F$, a quantifier-free Boolean formula defining a semialgebraic set $S \subset \R^n$

{\bf Output: $\mathcal{}$} $\mathcal{D}$, a CAD of $\R^n$ compatible with $S$ and satisfying the frontier condition\\

First, apply algorithm from Proposition~\ref{prop:collins-set} to $F$. We obtain a CAD $\mathcal E$ of $\R^n$
compatible with $S$.

The algorithm computes a sequence of pairs $({\mathcal E}_M,\ X_M)$ where $M \in \{0,1\}^{n-1}$ is an index, $\mathcal{E}_M$ is a CAD of $\R^n$, which is a refinement of $\mathcal{E}$, and $X_M \subset \R^n$ is a semialgebraic set contained in the union of all cells of ${\mathcal E}_M$ with indices $\lexeq (i_1,\ldots,i_{n-1},1)$. This family of cells will be denoted by $\mathcal{U}_M$. Note that index $M$ refers to the pair rather than its elements and $\mathcal{E}_M$ and $X_M$ are written with indices for convenience.

The sequence of pairs is computed recursively, starting with index $M=(1,1, \ldots ,1)$,
in descending order of indices $M=(i_1, \ldots, i_{n-1})$ with respect to $\lex$.
Let the initial pair $({\mathcal E}_{(1, \ldots, 1)},\ X_{(1, \ldots, 1)}) = ({\mathcal E},\ \emptyset)$. The algorithm will construct a sequence (from right to left):
\begin{equation}\label{eq:sequence}
({\mathcal E}_{\mathbf{0}}, X_{\mathbf{0}}) \lex \cdots \lex
({\mathcal E}_N, X_N) \lex ({\mathcal E}_M, X_M)
\lex \cdots
\lex ({\mathcal E} , \emptyset),
\end{equation}
where index $(0,0,\ldots,0)$ is abbreviated to $\mathbf{0}$ and $X_{\mathbf{0}}= \bigcup \{ C \in{\mathcal U}_{(0, \ldots ,0)} \}$.

For a given index $M=(i_1,\ldots,i_{n-1})$, assume that the algorithm has computed the pair $(\mathcal{E}_M,\ X_M)$.
Now we describe how the next pair, $({\mathcal E}_N,\ X_N)$, where $N=(j_1, \ldots ,j_{n-1})$ is the index immediately prior to $M$ with respect to $\lex$, is computed.
Applying algorithm from Proposition~\ref{prop:collins-sets} to ${\mathcal E}_M$ and $X_M$, we get a CAD
${\mathcal E}_{N}$ of $\R^n$ compatible with $X_M$ and every cell of ${\mathcal E}_M$ (see construction of $\mathcal{E}_{(0,1)}$ in step~\ref{step:one-zero} of the worked example).
Consider the family ${\mathcal A}_0$ of cells in ${\mathcal E}_{N}$ with indices $(i_1,\ldots,i_{n-1},0)$
(section cells)
and the family ${\mathcal A}_1$ of cells in ${\mathcal E}_{N}$ with indices $(i_1,\ldots,i_{n-1},1)$
(sector cells).
Using algorithm from Lemma~\ref{le:frontier}, compute the set
\begin{equation}\label{eq:front}
X_{N}:= \bigcup_{C \in {\mathcal A}_0}\fr (C) \cup \bigcup_{C \in {\mathcal A}_1}(\fr (C) \setminus (C_T \cup C_B)).
\end{equation}
Observe that for every cell $C \in {\mathcal A}_1$, $C_T$ and $C_B$ are section cells with index $(i_1,\ldots,i_{n-1},0)$ in ${\mathcal A}_0$.
According to Lemma~\ref{lem:fr-lex-less}, $X_{N}$ is contained in the union of cells in ${\mathcal E}_{N}$
with indices $\lex (i_1,\ldots,i_{n-1},0)$ (see construction of $X_{(1,0)}$ in step~\ref{step:one-one} of the worked example).

When the algorithm reaches the final pair $({\mathcal E}_{\mathbf{0}}, X_{\mathbf{0}})$ in the sequence,
it computes a CAD $\mathcal D$, using algorithm from Proposition~\ref{prop:collins-sets}, compatible
with $X_{\mathbf{0}}$ and every cell of ${\mathcal E}_{\mathbf{0}}$. The algorithm terminates by returning $\mathcal D$.
\medskip

\noindent{\bf Correctness.}\\
Let $L=(k_1, \ldots ,k_{n-1})$ be the index immediately prior
to $N=(j_1, \ldots ,j_{n-1} )$ with respect to $\lex$.
The algorithm computes ${\mathcal E}_L$ as the refinement of ${\mathcal E}_N$ compatible with $X_N$ (see Equation~\ref{eq:front}). Suppose that the algorithm has computed the final decomposition $\mathcal{D}$, as the refinement of $\mathcal{E}_{\mathbf{0}}$ compatible with $X_{\mathbf{0}}$ in the sequence \ref{eq:sequence}. Now construct a new sequence of decompositions (from left to right)
\begin{equation}
\label{eq:seq-ind-hyp}
\mathcal{E}'_{\mathbf{0}} \lex \cdots \lex \mathcal{E}'_L \lex \mathcal{E}'_N \lex \mathcal{E}'_M \lex \cdots \lex \mathcal{E}'_{(1,\ldots,1)}
\end{equation}
where the initial element $\mathcal{E}'_{\mathbf{0}}$ is the refinement of $\mathcal{E}_{\mathbf{0}}$ compatible with $\mathcal{D}$ and each $\mathcal{E}'_I$ is the refinement of $\mathcal{E}_I$ compatible with all cells in $\mathcal{E}'_J$ where $J$ is the index immediately prior to $I$ with respect to $\lex$.
We want to prove that every cell in $\mathcal{E}'_L$ with index $\lexeq (i_1,\ldots,i_{n-1},1)$ satisfies the frontier condition.

If $C$ is a cell in ${\mathcal E}_N$ with index $(i_1, \ldots ,i_{n-1} ,0)$ (section cell),
then $C$ is a union of cells of ${\mathcal E}_L$ and, hence, of $\mathcal{E}'_L$, with indices $\lexeq (i_1, \ldots ,i_{n-1} ,0)$.
Let $C' \subset C$ be one of the cells in this refinement of $C$ with index $(i_1, \ldots ,i_{n-1} ,0)$ and $B$ be
the union of cells contained in the refinement of $C$ with indices $\lex (i_1, \ldots ,i_{n-1} ,0)$.
According to Lemma~\ref{lem:section-split}, $\fr_C (C \setminus B)=B$ where $\fr_X(Y)$ denotes the frontier of $Y$ in $X$.
Hence, $\fr_C (C') \subset B$.
On the other hand, $\fr (C') \setminus \fr_C (C')$ is a subset of $\fr (C)$, which is a union of cells of $\mathcal{E}_L$ (and of the refinement $\mathcal{E}'_L$) with index $\lex (i_1,\ldots,i_{n-1},0)$.
By the induction hypothesis, all cells of $\mathcal{E}'_L$ with index $\lex (i_1,\ldots,i_{n-1},0)$ satisfy the frontier condition. In particular, $\fr(C)$ and $\fr(B)\subset \fr(C)$ is a union of cells of $\mathcal{E}'_L$.
It follows from the cylindrical structure of $\mathcal{E}'_L$ that $C'$ satisfies the frontier condition.

If $C$ is a cell in ${\mathcal E}_N$ with index $(i_1, \ldots ,i_{n-1} ,1)$ (sector cell),
then $C$ is a union of cells of ${\mathcal E}_L$ with indices $\lexeq (i_1, \ldots ,i_{n-1} ,1)$.
A similar argument to that used for section cells shows that each cell $C'$ in this union,
having index $(i_1, \ldots ,i_{n-1} ,1)$, satisfies the frontier condition in $\mathcal{E}'_L$. Note that when using
Lemma~\ref{lem:section-split} we consider sector cell $C$ as a graph of a constant function over itself.

Finally, each decomposition $\mathcal{E}'_I$, $I \in \{0,1\}^{n-1}$, in the sequence \ref{eq:seq-ind-hyp} coincides with $\mathcal{D}$. Indeed, $\mathcal{D}$ is a refinement of $\mathcal{E}_I$ and $\mathcal{E}'_I$ is a refinement of $\mathcal{E}_I$ compatible with all cells of $\mathcal{D}$. In other words, no refinement of $\mathcal{D}$ is required to obtain $\mathcal{E}'_{(1,\ldots,1)}$.
It follows that every cell of $\mathcal{D}$ satisfies the frontier condition.
\medskip

\noindent{\bf Complexity.}\\
The number of different indices $(i_1, \ldots, i_n)$, where each $i_k \in \{ 0,1 \}$, is $2^n$.
Therefore, the algorithm makes $O(2^n)$ ``steps'': computing successive pairs $({\mathcal E}_M, X_M)$ in the sequence (\ref{eq:sequence}).
On each step, passing from a pair $(\mathcal{E}_M,\ X_M)$ to the next pair $(\mathcal{E}_N,\ X_N)$, the algorithm applies Proposition~\ref{prop:collins} to $({\mathcal E}_M, X_M)$ and obtains a CAD
${\mathcal E}_N$ which is a refinement of ${\mathcal E}_M$ compatible with $X_M$.
Then Lemma~\ref{le:frontier} is applied to each cell with index $(i_1,\ldots,i_{n-1},i_n)$ in ${\mathcal E}_N$
to obtain $X_N$.

Suppose that there are $s_M$ polynomials defining ${\mathcal E}_M$ and $X_M$, having degrees at most $d_M$.
Then, according to Proposition~\ref{prop:collins-sets}, the CAD ${\mathcal E}_N$ is defined by
$s_N:=(s_Md_M)^{O(1)^n}$ polynomials of degrees $d_N=(s_Md_M)^{O(1)^n}$.
The number of cells in ${\mathcal E}_N$ is asymptotically the same: $(s_Md_M)^{O(1)^n}$.
Then the algorithm applies Lemma~\ref{lem:fr-lex-less} to compute $X_N$.
The complexity of this application is $(s_Md_M)^{n^2O(1)^n}$, which is asymptotically the same as $(s_Md_M)^{O(1)^n}$.
Thus, the overall complexity of this ``step'' is again $(s_Md_M)^{O(1)^n}$. The overall complexity is obtained by iterating this process for each index.
Given that $s_{(1, \ldots,1)}=s$ and $d_{(1, \ldots,1)}=d$, we conclude that
the complexity of the algorithm is
$$
(sd)^{O(1)^{n2^n}}.
$$
This is also an upper bound on the number of cells in $\mathcal D={\mathcal E}_{\mathbf{0}}$, the number of polynomials defining cells and their degrees.
\end{proof}

\subsection{An Aside on Constructing the Intermediate Decompositions}

Throughout this algorithm we construct ``a CAD compatible with $F$ and all cells of $\mathcal{E}$''.
This can be achieved by constructing a CAD, using Proposition~\ref{prop:collins}, such that each cell has constant sign on the polynomials defining $F$ and all cells of $\mathcal{E}$. Such a CAD is clearly compatible with $F$ and all cells of $\mathcal{E}$. However, it may include some cells, outside $\cl(S)$, which we are not interested in.

An alternative is to use the algorithm described in Proposition~\ref{prop:collins-sets} to construct a CAD compatible with $\{F, C_1,\ldots,C_r\}$ where $C_i,1\le i \le r$ is a cell of $\mathcal{E}$ such that $\bigcup_{1\le i \le r} C_i = \R^n$. Observe that this CAD is compatible with the required sets.

Both options have the same asymptotic complexity. Therefore, the choice of CAD subroutine does not change the complexity bound obtained in the proof of Theorem~\ref{th:main}. In the proof, Proposition~\ref{prop:collins-set}, constructing a CAD compatible with a semialgebraic set $S\subset \R^n$, is used to construct the initial CAD and Proposition~\ref{prop:collins-sets}, constructing a CAD compatible with a family of sets, is used for constructing the refinements of intermediate decompositions.

It is assumed that every cell in $\mathcal{E}$ can be represented as a quantifier-free Boolean formula. This is always possible by Thom's Lemma, but, as shown by Brown in \cite{brown99}, the polynomials produced during the construction of $\mathcal{E}$ may not be sufficient to do so. It is possible to obtain the required polynomials by using Collins' extended projection. Brown \cite{brown99} also presents an algorithm that ensures every cell can be represented by a formula containing projection polynomials and, possibly, a small number of derivatives. This may be the preferred option since fewer extra polynomials are needed.

\section{Worked Example}
\label{sec:worked-example}

We demonstrate the algorithm from Theorem~\ref{th:main} by applying it to the set
$$
S := \{ -1 < x < 1, \vert x\vert <y <1, z=\vert x/y\vert \}
$$
from Example~\ref{example:top-bottom-not-cylindrical}. Construct a CAD $\mathcal{E}$ of $\R^3$ compatible with $S$ (see Table~\ref{table:cells-of-e}). Observe that $S$ is a single 2-dimensional cylindrical $(1,1,0)$-cell of $\mathcal{E}$.
\begin{enumerate}
\item \label{step:one-one}
The initial pair $(\mathcal{E}_{(1,1)}, X_{(1,1)}) = (\mathcal{E}, \emptyset)$. We now describe how the next pair in the sequence (\ref{eq:sequence}) is computed.
$\mathcal{E}_{(1,0)}$, the CAD compatible with all cells of $\mathcal{E}_{(1,1)}= \mathcal{E}$ and $X_{(1,1)} = \emptyset$, is equal to $\mathcal{E}$. $X_{(1,0)}$ is contained in the union of cells of $\mathcal{E}_{(1,0)}$ with index $\lex (1,1,0)$ and includes
$$
\fr(S) = S_B \cup S_T \cup \{ (-1,1,1), (1,1,1) \}
$$
where
$$
S_T = \{ -1 < x < 1, y = 1, z = \vert x \vert \}
$$
and
\begin{align*}
S_B = &\{ -1 < x < 0, y = -x, z = 1 \} \cup \\
& \{ x = y = 0, 0 \le z \le 1 \} \cup\\
& \{ 0 < x < 1, y = x, z = 1 \}.
\end{align*}
Observe that $C_T$ could be a $(1,0,0)$-cell (although it is not a cell in $\mathcal{E}_{(1,0)}$), but $C_B$ cannot be a cell in any cylindrical decomposition since it is not the graph of a continuous function.

\item \label{step:one-zero}
Given $(\mathcal{E}_{(1,0)}, X_{(1,0)})$, compute the next pair.
$\mathcal{E}_{(0,1)}$ is compatible with every cell of $\mathcal{E}_{(1,0)} = \mathcal{E}$ and $X_{(1,0)} \supset \fr(S)$.
The blow-up subset $\{ x = y = 0, 0 \le z \le 1 \}$ in $X_{(1,0)}$ results in a refinement of the decomposition induced by $\mathcal{E}_{(1,0)}$ on $\R^1$ such that it includes the cells $\{-1 < x < 0\}, (0), \{ 0 < x < 1 \}$. Thus $\mathcal{E}_{(0,1)}$ includes a cell $S'' = \{ x = 0, 0 < y < 1, z = 0 \}$ and $S$ is split into three cells $S', S'', S'''$ with indices $(1,1,0),(1,0,0)$ and $(1,1,0)$ respectively.
Cells $C_1,C_2,C_3$ and $C_4$ of $\mathcal{E}_{(1,0)} = \mathcal{E}$ are also refined in this step so that they are compatible with $\fr(S)$ (see Table~\ref{table:cells-of-e-refined}). As argued in the Correctness section, frontiers $\fr_S(S')$ and $\fr_S(S''')$ coincide with $S''$, so no further refinements of cells with index $(1,1,1)$ and $(1,1,0)$ are needed.
$X_{(0,1)}$ contains $\fr(C'_2), \fr(C'''_2), \fr(C'_3)$ and $\fr(C'''_3)$.

\item \label{step:zero-one}
Observe that, in this particular case, $X_{(0,1)}$ is already a union of cells of $\mathcal{E}_{(0,1)}$, so $\mathcal{E}_{(0,0)} = \mathcal{E}_{(0,1)}$ as no refinement is needed.
$X_{(0,0)}$ contains $\fr(S')$.

\item $X_{(0,0)}$ is already a union of cells, so no refinement of $\mathcal{E}_{(0,0)} = \mathcal{E}_{(0,1)}$ is needed. The algorithm terminates and returns $\mathcal{D} = \mathcal{E}_{(0,1)}$.
\end{enumerate}

\begin{table}
\caption{All cells of the CAD $\mathcal{E}$, the initial CAD computed in Section~\ref{sec:worked-example}.}

\begin{tabular}{c | c | c}\label{table:cells-of-e}
Label & Index & Formula\\
\hline
& $(1,1,1)$ & $\{ x < -1 \}$\\
$C_1$& $(0,1,1)$ & $\{ x = -1 \}$\\
& $(1,1,1)$ & $\{ -1 < x < 1, y < \vert x \vert \}$\\
$C_2$& $(1,0,1)$ & $\{ -1 < x < 1, y = \vert x \vert \}$\\
& $(1,1,1)$ & $\{ -1 < x < 1, \vert x \vert < y < 1, z < \vert x/y \vert \}$\\
$S$ & $(1,1,0)$ & $\{ -1 < x < 1, \vert x \vert < y < 1, z = \vert x/y \vert \}$\\
& $(1,1,1)$ & $\{ -1 < x < 1, \vert x \vert < y < 1, z > \vert x/y \vert \}$\\
$C_3$ & $(1,0,1)$ & $\{ -1 < x < 1, y = 1 \}$\\
& $(1,1,1)$ & $\{ -1 < x < 1, y > 1 \}$\\
$C_4$& $(0,1,1)$ & $\{ x = 1 \}$\\
& $(1,1,1)$ & $\{ x > 1 \}$\\
\end{tabular}
\end{table}

\begin{table}
\caption{Cells of the CAD $\mathcal{E}_{(0,1)}$ computed in step \ref{step:one-zero} of Section~\ref{sec:worked-example}. Note that only cells which are part of $\cl(S)$ are listed.}

\begin{tabular}{c | c | c}\label{table:cells-of-e-refined}
Label & Index & Formula\\
\hline
${C'_1}$ & $(0,0,0)$ & $\{ x = 0, y = 1, z = 1 \}$\\
${C'_2}$ & $(1,0,0)$ & $\{ -1 < x < 0, y = -x, z = 1 \}$\\
${C''_{2,1}}$ & $(0,0,0)$ & $\{ x = y = 0, z = 0 \}$\\
${C''_{2,2}}$ & $(0,0,1)$ & $\{ x = y = 0, 0 < z < 1 \}$\\
${C''_{2,3}}$ & $(0,0,0)$ & $\{ x = y = 0, z = 1 \}$\\
${C'''_2}$ & $(1,0,0)$ & $\{ 0 < x < 1, y = x, z = 1 \}$\\
$S'$ & $(1,1,0)$ & $\{ -1 < x < 0, -x < y < 1, z =  -x/y \}$\\
$S''$ & $(0,1,0)$ & $\{ x = 0, 0 < y < 1, z = 0 \}$\\
$S'''$ & $(1,1,0)$ & $\{ 0 < x < 1, x < y < 1, z = x/y \}$\\
$C'_3$ & $(1,0,0)$ & $\{ -1 < x < 0, y = 1, z = -x \}$\\
$C''_3$ & $(0,0,0)$ & $\{ x = 0, y = 1, z = 0\}$\\
$C'''_3$ & $(1,0,0)$ & $\{ 0 < x < 1, y = 1, z = x \}$\\
${C'_4}$ & $(0,0,0)$ & $\{ x = 1, y = 1, z = 1 \}$\\
\end{tabular}
\end{table}

Well-known implementations of cylindrical algebraic decomposition such as Brown's {\em QEPCAD B} \cite{brownQepcad} and Maple's {\em CylindricalAlgebraicDecompose} \cite{chen2014} are unable to obtain the frontier condition in general. In particular, when applied to set $S$ from Example~\ref{example:top-bottom-not-cylindrical}, both {\em CylindricalAlgebraicDecompose} and {\em QEPCAD B} create three cells above the origin: $B' = \{ x = y = 0, z < 0\}$, $B'' = \{ x = y = 0, z = 0\}$ and $B''' = \{ x = y = 0, z > 0\}$. Observe that $B''' \cap \fr(S) \ne \emptyset$, but $B''' \not \subset \fr(S)$ so the frontier of $S$ is not a union of cells in the decomposition.

\section{Pseudo-code}

We now present the algorithm from Theorem~\ref{th:main} as pseudo-code. First define the following basic subroutines.

\begin{itemize}
\item
${\mathcal B}:= {\mathcal A}\ \&\ G$ takes a CAD $\mathcal A$ of $\R^n$ and a semialgebraic set $G \subset \R^n$,
and returns a CAD $\mathcal B$ of $\R^n$ compatible with $G$ and all cells of $\mathcal A$.
This subroutine follows from the algorithm in Proposition~\ref{prop:collins-sets}.

\item $Y := \fr(X)$ takes a set $X \subset \R^n$ and returns its frontier.
This subroutine follows from the algorithm in Lemma~\ref{le:frontier}.

\item $N := \rm{decrement}(M)$ takes an index $M \in \{0,1\}^k$ and returns the index $N \in \{0,1\}^k$ immetiately prior to $M$ with respect to $\lex$.
\end{itemize}

\begin{algorithm}

${\mathcal D} := \Frontier(F)$

\textbf{Input:}
$F$: a quantifier-free Boolean formula representing a semialgebraic set $S \subset \R^n$.

\textbf{Output:}
CAD $\mathcal D$ of $\R^n$ compatible with $S$ and satisfying the frontier condition.
\medskip

\begin{itemize}
\item
{\tt let} ${\mathcal E}:= \R^n\ \&\ S$, where $\R^n$ is a trivial CAD of $\R^n$,

\item $M := (1,\ldots,1) \in \{0,1\}^{n-1}$,

\item
$({\mathcal E}_M, F_M):=({\mathcal E}, \emptyset)$.

\item
{\tt while} $(0,\ldots,0) < M$

\noindent {\tt do}

\begin{itemize}
\item {\tt let} $N := {\rm decrement(M)}$.

\item
Construct $X_N$ from cells of $\mathcal E_M$ using ``$\fr$'' subroutine and formula~\ref{eq:front}.

\item
{\tt let} ${\mathcal E}_N := {\mathcal E}_M\ \&\ X_N$,

\item
$M := N$.
\end{itemize}

\noindent {\tt end while}

\item {\tt return} ${\mathcal E}_M$

\end{itemize}
\end{algorithm}

\section{Generalisations and Further Work}

\subsection{Pfaffian Functions}
\label{sec:pfaffian}

Theorem~\ref{th:main} can be extended to restricted sub-Pfaffian sets as described by Gabrielov and Vorobjov in \cite{gv04}.
In particular, it can be used to prove the existence of decompositions with frontier condition compatible with semialgebraic sets defined by {\em fewnomials} \cite[Section 2.6]{gv04}, whose structure is destroyed by a change of coordinates, and obtain a bound on the number of cells in these decompositions.

The idea of the proof is still valid. Only a modification to the subroutines for computing the frontier, and for constructing a classical CAD, is needed. For the former, replace the quantifier elimination algorithm in Lemma~\ref{le:frontier} with the result of \cite[Section 5]{gv04}. In the latter case, the classical CAD algorithm from Proposition~\ref{prop:collins} should be replaced by the main result of \cite{gv01} (see also \cite[Section 7]{gv04}). The rest of the proof can be reproduced almost identically.

Note that it is not currently clear how this result can be implemented since in the algorithm for cylindrical decomposition in \cite{gv01} the oracle is needed to decide whether a sub-Pfaffian set is empty or not.

\subsection{Further Work}

The complexity upper bound of the algorithm from Theorem~\ref{th:main} is significantly worse than the bound
for classical CAD algorithm in \cite{collins1975} and \cite{wuthrich2005}.
This is caused by the parameter of its recursive loop: the index of a cylindrical cell, which is exponential in the ambient dimension.
It is difficult to see another parameter that could make the recursion significantly shorter. Thus, if any progress is to be made towards a better asymptotic complexity, the method used may need to be based on a completely different ideas.
On the other hand, there is a lower bound on complexity of the classical CAD algorithm due to Davenport and Heintz \cite{davenportHeintz1988}. However, it is not yet known whether the lower bound for CAD with frontier condition is greater than that of the classical CAD.
Another strand of research could be to explore this question, possibly by attempting to raise this lower bound for CAD with frontier condition.

Another improvement in the algorithm might come from a different subroutine for computing the frontier
of a cylindrical cell.
The subroutine from Lemma~\ref{le:frontier} works for any semialgebraic set and does not take advantage of its
cylindrical structure.
We may be able to exploit this structure by factorising the equational part of the formula
representing a cell into irreducible components, generalising a method in \cite{lazard10}.

Finally, since some initial CAD is refined repeatedly, the classical CAD algorithm from Proposition~\ref{prop:collins} could be replaced with an incremental algorithm, E.g., the algorithm presented by Kremer and \'Abrah\'am in \cite{kremer2020}. It is not clear whether substituting the CAD algorithm will reduce the complexity bound. The set $X_N$ is still computed at each step, using a singly exponential algorithm, and the initial CAD may need to be refined at every step and this refinement may result in the incremental algorithm backtracking all the way to the decomposition induced on $\R^1$. However, this change would be very useful in practice as unnecesary CAD recomputations could be avoided, E.g., if $\mathcal{E}_N$ is already compatible with the set $X_N$, the algorithm in its current form will compute the refinement even though it is not needed.

\subsubsection{Acknowledgements}

Thanks to Prof. Nicolai Vorobjov, my PhD supervisor, for all his help.
Thanks to Prof. James Davenport, Dr Matthew England and Tereso Del Río Almajano for insightful informal discussions.

%
%
%
\bibliographystyle{splncs04}
\bibliography{references}
\end{document}